\documentclass[12pt]{article}
\usepackage{amssymb,amsmath,latexsym}
\usepackage{graphicx}

\newtheorem{theorem}{Theorem}
\newtheorem{definition}{Definition}

\newtheorem{lemma}{Lemma}

\begin{document}
\title{Geometric flow of $G_{2}$-structures on $C(S^3\times S^3)$\footnote{This research was partially supported
by Grant of the Russian Federation for the State Support of
Researches (Agreement No 14.B25.31.0029).}}
\author{Khazhgali Kozhasov}
\date{July 2014}
\maketitle
\begin{abstract}
\begin{center}
    {We introduce a first order flow of $G_{2}$-structures and  construct its explicit
    solution in case of a cone over $S^3\!\times\! S^3$. Also we prove for this situation
    that starting from certain initial datum the flow deforms corresponding to
    $G_2$-structure metric to a conic metric up to homotheties.}
\end{center}
\end{abstract}

\section{Introduction}

Theory of flows of $G_2$-structures has appeared quite recently. The
most famous examples of these flows are Laplacian flow that was
introduced by R. Bryant in \cite{3} and S. Karigiannis' General flow
\cite{4}. It's not clear whether there exists a 'distinguished' flow
of $G_2$-structures that would lead to a parallel structure on
manifolds satisfying some (still unknown) conditions.

In this article we construct a first order flow of $G_2$-structures
that possesses interesting solutions in special cases. We find
explicit solution of this flow in case when our manifold is a cone
over $S^3\!\times\! S^3$ given by certain family of
$G_2$-structures.

For mentioned manifold we show that corresponding to $G_2$-structure
metric satisfying some conditions evolves along the flow to a conic
metric up to rescalings at every moment of time.

The author is grateful to Yaroslav Bazaikin for useful discussions
and advices.

\newpage
\section{Definition of the geometric flow}

\begin{figure}[h!]
\center{\includegraphics[width=0.5\linewidth]{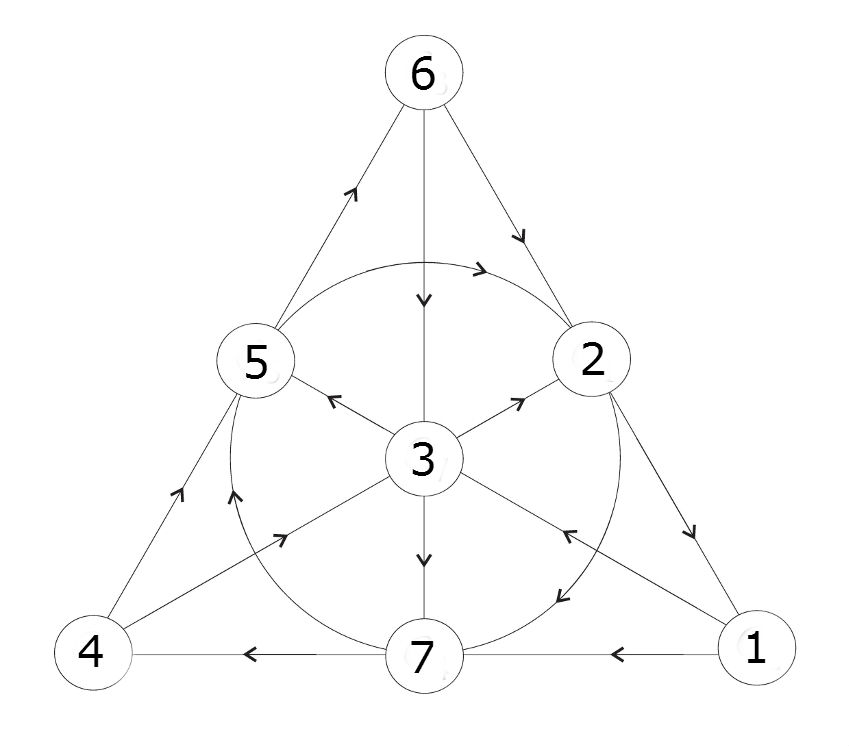}}\\
\caption{Fano plane} \label{ris:fano}
\end{figure}

Let's consider 8-dimensional octonion algebra $\mathbb{O}$ with a
basis $1$, $e_{1}$, $e_{2}$, ... , $e_{7}$ and a multiplication law
shown on fig.~\ref{ris:fano}. We will identify imaginary
Im$(\mathbb{O})$=$Span(e_{1}$, $e_{2}$, ... , $e_{7})$ octonions and
$\mathbb{R}^7$. Multiplication $\circ$ in $\mathbb{O}$ defines
positive definite $<\cdot,\cdot>$ and cross $\cdot\times\cdot$
products on pairs of vectors $u$, $v$ $\in$ $\mathbb{R}^7\cong
Im(\mathbb{O})$ as follows: $$<u,v>\ =-Re(u\circ v)$$
$$u\times v=Im(u\circ v)$$ These products allow us to define a 3-form (associative form) $\phi$ on $\mathbb{R}^7$
by the formula
$$\phi(u,v,w)=<u\times v,w>$$
\\This form is non-degenerate in the next sense: for any $x$,
$y$ $\in\mathbb{R}^7$
$$(x\lrcorner\phi)\wedge (y\lrcorner\phi)\wedge\phi\neq0$$
Written in the basis $e_{1}$, $e_{2}$, ... , $e_{7}$ it looks as
follows:
$$\phi=e^{456}+e^{621}+e^{174}+e^{527}+e^{637}+e^{135}+e^{432}$$
where by $e^{ijk}$ we denote the basic form $e^i\wedge e^j\wedge
e^k$ and $e^i(e_{j})=\delta^i_{j}$.
\begin{definition}
The subgroup of $GL_7$ that fixes (with respect to canonical action
$GL_{7}\hookrightarrow\wedge^3(\mathbb{R}^7)$) the non-degenerate
form $\phi$ is called the group $G_2$.
\end{definition}
It is 14-dimensional simple Lie group of type $G_2$. The orbit
$\wedge^3_{+}(\mathbb{R}^7)$ of $\phi$ that consists of
non-degenerate 3-forms is open in $\wedge^3(\mathbb{R}^7)$ because
$$dim\ \Lambda^3_{+}(\mathbb{R}^7)=dim\ GL_{7}-dim\
G_{2}=dim\ \Lambda^3(\mathbb{R}^7)$$

Henceforward $M$ is a 7-dimensional manifold. Differential 3-form
$\varphi$ on $M$ is called non-degenerate if $\varphi(x)$ is
non-degenerate as a 3-form on $T_{x}(M)$ $\forall x\in M$.

\begin{definition}
Any non-degenerate 3-form $\varphi$ on $M$ will be called a
$G_{2}$-structure on $M$.
\end{definition}
In a neighborhood of a point $p\in M$ such local coordinates could
be chosen that in these coordinates $\varphi(p)$ coincides with
above described associative form $\phi$. The set $\Lambda^3_{+}(M)$
of non-degenerate 3-forms on $M$ is just the set of smooth sections
of the bundle $\Lambda^3_{+}(TM)$ over $M$ with a fiber
$\Lambda^3_{+}(\mathbb{R}^7).$ It is known that $M$ admits a
$G_2$-structure if and only if the first two Stiefel-Whitney classes
of $M$ vanish.

$G_{2}$-structure on $M$ allows one to define riemannian metric
$g=g_{\varphi}$ on $M$ as it follows. In local coordinates $x^1$,
$x^2$, ... , $x^7$ in a neighborhood of point $p\in M$ let's define
the tensor field $B$ by the rule
$$B_{ij}\ dx^1 \wedge ... \wedge dx^7 = \frac{\partial}{\partial
x^i}\lrcorner\varphi\wedge\frac{\partial}{\partial
x^j}\lrcorner\varphi\wedge\varphi.$$ Then the metric $g$ is defined
by the formula
$$g_{ij}=\frac{1}{6^\frac{2}{9}}\frac{B_{ij}}{det(B)^\frac{1}{9}}.$$
If one choose local coordinates such that in these coordinates
$\varphi(p)$ will be expressed as associative form $\phi$ the metric
will have an Euclidean form $g_{ij}(p)=\delta_{ij}$.
\begin{definition}
If for the $G_2$-structure $\nabla\varphi=0$, where $\nabla$ is a
Levi-Civita connection of the metric $g=g_{\varphi}$, then
($M$,$\varphi$) is called a $G_2$-manifold.
\end{definition}
\textbf{Remark}:\textit{The equation $\nabla\varphi=0$ is a very
non-linear and thus it's so difficult to find on $M$ a parallel
$G_{2}$-structure.}\\ Any $G_2$-manifold has its holonomy group
contained in $G_{2}$. It's interesting that $G_2$-manifolds are
always Ricci-flat:
\begin{theorem}[Bonan, \cite{1}]
$$\nabla\varphi=0\Rightarrow Ric(g_{\varphi})=0$$
\end{theorem}
There is a useful criteria for $G_2$-structure to be parallel:
\begin{theorem}[Fernandez and Gray, \cite{2}]
$$\nabla\varphi=0\Leftrightarrow
d\varphi=0\ \textrm{and}\ \delta\varphi=0$$ where $\delta$ is the
conjugate operator to the De Rham differential $d$ w.r.t. the metric
$g_{\varphi}$.
\end{theorem}

\begin{definition}
Let $\varphi=\varphi(t)$ be a smooth family of
$G_{2}$-strcutures on $M$. The flow of $G_{2}$-structures is a
system of evolutionary differential equations for components of
$\varphi$ in a basis $dx^i\wedge dx^j\wedge dx^k$:
$$\frac{\partial\varphi_{ijk}}{\partial t}=F(\varphi)_{ijk}$$where
$F(\varphi)_{ijk}$ --- some, generally speaking, differential (in
sense of only space variables) expressions on components of
$\varphi$.
\end{definition}
Examples of flows of $G_2$-structures are\\
--- the Laplacian flow \cite{1} that was introduced by R. Bryant
$$\frac{\partial\varphi}{\partial t}=\Delta_{g}\varphi$$ where
$\Delta_{g}=d\delta+\delta d$ is a Hodge Laplacian of the metric $g_{\varphi}$,\\
--- the General flow \cite{3} introduced by S. Karigiannis
$$\frac{\partial\varphi_{ijk}}{\partial
t}=h^{l}_{i}\varphi_{ljk}+h^l_{j}\varphi_{ilk}+h^l_{k}\varphi_{ijl}+X^l(*\varphi)_{lijk}$$
where $h_{ij}$ is a symmetric tensor on $M$, $X^k$ is a vector field
on $M$, $*$ is a Hodge star operator of the metric $g_{\varphi}$ .
\begin{definition}
The geometric flow of $G_2$-structures on $M$ is the equation
\begin{equation*}
    \frac{\partial\varphi}{\partial t}\wedge X=d\varphi
\end{equation*}
where $X$ is a differential 1-form on $M$ that does not depend on
time.
\end{definition}

\section{Geometric flow on $C(S^3\!\times\! S^3)$}\label{section}

Let's consider a cone over $S^3\!\times\! S^3$. There are a frame of
left-invariant vector fields on the Lie group $S^3=SU(2)$ such that
at the unit of $SU(2)$ it looks as follows
$$\xi^1=\begin{pmatrix} i & 0 \\ 0 & -i \end{pmatrix},\
\xi^2=\begin{pmatrix} 0 & 1 \\ -1 & 0
\end{pmatrix},\ \xi^3=\begin{pmatrix} 0 & i\  \\ \ i & 0
\end{pmatrix}.$$ The Lie algebra generated by these vector fields has the multiplication law
$$[\xi^1,\xi^2]=2\xi^3,\ [\xi^2,\xi^3]=2\xi^1,
\ [\xi^3,\xi^1]=2\xi^2.$$

If $\eta_1$, $\eta_2$, $\eta_3$ are dual to $\xi^1$, $\xi^2$,
$\xi^3$ 1-forms then, by Cartan formula,
$$d\eta_1=-2\eta_2\wedge\eta_3,\ d\eta_2=-2\eta_3\wedge\eta_1,\
d\eta_3=-2\eta_1\wedge\eta_2.$$

Let $\eta_1$, $\eta_2$, $\eta_3$, $\tilde{\eta}_1$ and
$\tilde{\eta}_2$, $\tilde{\eta}_3$ be left-invariant coframes on the
first and on the second sphere of the Cartesian product
$S^3\!\times\! S^3$ correspondingly, let $dr$ be a standart 1-form
on $\mathbb{R}$. Now define the next 1-forms on $C(S^3\times S^3)$
$$e^1=A(r)(\eta_1+\tilde{\eta_1}),$$
$$e^2=A(r)(\eta_2+\tilde{\eta_2}),$$
$$e^3=A(r)(\eta_3+\tilde{\eta_3}),$$
$$e^4=B(r)(\eta_4-\tilde{\eta_4}),$$
$$e^5=B(r)(\eta_5-\tilde{\eta_5}),$$
$$e^6=B(r)(\eta_6-\tilde{\eta_6}),$$
$$e^7=dr$$
where $A(r)$ and $B(r)$ --- some strictly positive functions, $r>1$.\\
In the basis $e^i,\ i=1,..,7$ the $G_{2}$-structure and the
corresponding metric are given by the formulas
\begin{equation}\label{varphi}
    \varphi=e^{456}+e^{621}+e^{174}+e^{527}+e^{637}+e^{135}+e^{432},
\end{equation}
\begin{equation}\label{metric}
g=dr^2+\sum\limits_{i=1}^3A^2(\eta_{i}+\tilde{\eta_{i}})^2
+\sum\limits_{j=1}^3B^2(\eta_{j}-\tilde{\eta_{j}})^2
\end{equation}
Parallelness of just described $G_2$-structure was studied in \cite{5}.\\
\textbf{Remark}: \textit{There is a one-to-one correspondence
between the class of $G_{2}$-structures \ref{varphi} on
$C(S^3\!\times\! S^3)$  and the class of metrics \ref{metric} on
$C(S^3\times S^3)$.}

If $A$ and $B$ are functions smoothly depending on $t$, namely we
deal with a smooth family of $G_{2}$-structures
\begin{multline*}
$$\varphi=\varphi(t)=B^3(\eta_4-\tilde{\eta_4})\wedge(\eta_5-\tilde{\eta_5})\wedge
(\eta_6-\tilde{\eta_6})+A^2B(\eta_6-\tilde{\eta_6})\wedge(\eta_2+\tilde{\eta_2})\wedge
(\eta_1+\tilde{\eta_1})\\+AB(\eta_1+\tilde{\eta_1})\wedge
dr\wedge(\eta_4-\tilde{\eta_4})+
AB(\eta_5-\tilde{\eta_5})\wedge(\eta_2+\tilde{\eta_2})\wedge
dr+AB(\eta_6-\tilde{\eta_6})\wedge (\eta_3+\tilde{\eta_3})\wedge
dr\\+A^2B(\eta_1+\tilde{\eta_1})\wedge(\eta_3+\tilde{\eta_3})\wedge
(\eta_5-\tilde{\eta_5})+A^2B(\eta_4-\tilde{\eta_4})\wedge(\eta_3+\tilde{\eta_3})\wedge
(\eta_2+\tilde{\eta_2}),$$
\end{multline*}
then we have a
\begin{lemma}\label{lemma1} The geometric flow
\begin{equation}\label{flow}
    \frac{\partial\varphi}{\partial t}\wedge dr=d\varphi
\end{equation}
is equivallent to the system of PDE
\begin{equation}\label{system}
    \begin{cases}
        A=2BB_{x}
        \\
        8BB_{xx}+12B_{x}^2=1
    \end{cases}
\end{equation}
where $x=t+r$, $y=r-t$.
\end{lemma}
Proof of \ref{lemma1} is in Appendix.\\
\begin{lemma}
General solution of \ref{system} is given by the next expressions
\begin{equation}\label{general solution}
    \begin{cases}
        A=2B\sqrt{\frac{1}{12}+\frac{f(y)}{B^3}}
        \\
        x=\int\limits_{0}^{B}\frac{db}{\sqrt{\frac{1}{12}+\frac{f(y)}{b^3}}}+h(y)
    \end{cases}
\end{equation}
where $f(y)$ and $h(y)$ are arbitrary smooth functions.
\end{lemma}
\textit{Proof:} $B$ and  $B_{x}$ are strictly positive functions,
because any degeneration means that \ref{metric} is not a riemannian
metric. The equation
\begin{equation}\label{equation1}
    8BB_{xx}+12B_{x}^2=1,
\end{equation}
is equivallent to
$$(B_{x}^2B^3)_{x}=\frac{1}{12}(B^3)_{x}.$$
By integrating the last equation we get
\begin{equation}\label{equation2}
    B_{x}^2=\frac{1}{12}+\frac{f(y)}{B^3}
\end{equation}
where $f(y)$ --- arbitrary smooth function.\\
Let's fix $y\prime$ and construct solution of \ref{equation1} along
the characteristic $y=y\prime$. $B_{x}>0$ so we have
\begin{equation}\label{equation3}
    \frac{dx}{dB}=\frac{1}{\sqrt{\frac{1}{12}+\frac{f(y\prime)}{B^3}}}.
\end{equation}
and
\begin{equation*}
    x=\int\limits_{0}^{B}\frac{db}{\sqrt{\frac{1}{12}+\frac{f(y\prime)}{b^3}}}+h(y\prime).
\end{equation*}
 This expression is valid for any $y$ because $y\prime$ was an arbitrary value:
\begin{equation}\label{equation4}
    x=\int\limits_{0}^{B}\frac{db}{\sqrt{\frac{1}{12}+\frac{f(y)}{b^3}}}+h(y),
\end{equation}
Combining \ref{equation2} with the first equation of \ref{system} we
end the proof. $\Box$

\begin{definition}
Let $g=g(t,r)$ be a continuous family of metrics on
$C(S^3\times S^3)$. We say that $g$ converges to a metric
$g_{\infty}$ and write $g\rightarrow g_{\infty}$ as $t\rightarrow
+\infty$, if
\begin{equation}\label{definition of convergence}
    \forall\ K>1\sup_{1<r\leq K}\left|g(r,t)-g_{\infty}(r)\right|\rightarrow\ 0\ as\ t\rightarrow +\infty.
\end{equation}
\end{definition}
The main result of this paper is the next
\begin{theorem}\label{theorem}
For bounded $f(y)$ and $h(y)$ such that $f(y)\geq0$, $h(y)<y$ the
metric $g$ corresponding to solution \ref{general solution} of the
flow \ref{flow} satisfy the next condition: $\frac{g}{(t+1)^2}$
converges to a conical metric $ds^2+s^2\cdot g_{S^3\times S^3}$,
$g_{S^3\times S^3}$
--- metric on $S^3\times S^3$ that does not depend on $s$.
In other words: homothety class of the metric $g$ converges to a
homothety class of $g_{\infty}$ in the sense \ref{definition of
convergence}.
\end{theorem}
\textit{Proof:} As $r>1,\ t\geq 0$ then at each time $t$ $x>t+1$.
Further, instead of variables $x$ and $y$ we will sometimes use
$s=\frac{x}{t+1}$ and $t$, where $s>1$
--- a space variable of the limit metric.
Let's prove that
\begin{equation}\label{convergence}
    \sup_{s>1}|\frac{B(s,t)}{t+1}-\frac{s}{\sqrt{12}}|\rightarrow 0,
    \ \textrm{as}\ t\rightarrow +\infty
\end{equation}
Norming by $t+1$ is exactly homothety of corresponding metric.

By constituing $t=0$ into \ref{equation4} we get
$$r=\int\limits_{0}^{B|_{t=0}(r)}\frac{db}{\sqrt{\frac{1}{12}+\frac{f(r)}{b^3}}}+h(r),$$
$$h(y)=y-\int\limits_{0}^{B|_{t=0}(y)}\frac{db}{\sqrt{\frac{1}{12}+\frac{f(y)}{b^3}}}.$$
\textbf{Remark}: \textit{The condition $h(y)<y$ corresponds to
positiveness of $B|_{t=0}$}

By the mean value theorem for any $B>0$ there exists $B\prime$ such
that
$$x=\int\limits_{0}^{B}\frac{db}{\sqrt{\frac{1}{12}+\frac{f(y)}{b^3}}}+h(y)=
\frac{B}{\sqrt{\frac{1}{12}+\frac{f(y)}{B\prime^3}}}+h(y).$$ If
$x\rightarrow +\infty$ then $B,B\prime\rightarrow +\infty$ because
$f(y)\geq0,\ h(y)$ are bounded and
$B_{x}\rightarrow\frac{1}{\sqrt{12}}$ as $B\rightarrow +\infty$.

So, when $x\approx +\infty$
$B\approx\frac{x}{\sqrt{12}}-\frac{h(y)}{\sqrt{12}}$

Let's justify the convergence \ref{convergence}. Firstly, let's
notice that when $t$ is fixed we have
$$\frac{\partial}{\partial s}\left(\frac{B(s,t)}{t+1}-\frac{s}{\sqrt{12}}\right)=
\sqrt{\frac{1}{12}+\frac{f(y)}{B^3}}-\frac{1}{\sqrt{12}}\geq0$$ è
$$\sqrt{\frac{1}{12}+\frac{f(y)}{B^3}}-\frac{1}{\sqrt{12}}\rightarrow0,\
\textrm{as}\ s=\frac{x}{t+1}\rightarrow +\infty\Leftrightarrow
x\rightarrow +\infty$$ because $g$ is bounded and $\geq 0$.

Thus we get
\begin{equation}\label{estimate}
    \sup_{s>1}\left|\frac{B(s,t)}{t+1}-\frac{s}{\sqrt{12}}\right|=\frac{\left|h(y)\right|}{\sqrt{12}(t+1)}\rightarrow0\
    \textrm{as}\ t\rightarrow +\infty
\end{equation}

Now we can show that
\begin{equation}\label{convergence1}
    \forall\ K>1\ \sup_{1<s\leq K}\left|\frac{B(s,t)^2}{(t+1)^2}-\frac{s^2}{12}\right|\rightarrow 0,
    \ \textrm{as}\ t\rightarrow +\infty
\end{equation}
\begin{multline*}
    \sup_{1<s\leq K}\left|\frac{B(s,t)^2}{(t+1)^2}-\frac{s^2}{12}\right|\leq
    \sup_{1<s\leq K}\left|\frac{B(s,t)}{t+1}-\frac{s}{\sqrt{12}}\right|
    \left|\frac{B(s,t)}{t+1}+\frac{s}{\sqrt{12}}\right|\leq\\
    \frac{\left|h(y)\right|}{\sqrt{12}(t+1)}\left(\frac{\left|h(y)\right|}
    {\sqrt{12}(t+1)}+2\sup_{1<s\leq K}\frac{s}{\sqrt{12}}\right)\rightarrow0\
    \textrm{as}\ t\rightarrow +\infty
\end{multline*}

Finally, we prove the convergence \ref{definition of convergence} of
the metric
\begin{multline*}
    \frac{g}{(t+1)^2}=\frac{1}{(t+1)^2}\left(dx^2+\sum\limits_{i=1}^3A^2(\eta_{i}+\tilde{\eta_{i}})^2
    +\sum\limits_{j=1}^3B^2(\eta_{j}-\tilde{\eta_{j}})^2\right)=ds^2+
    \sum\limits_{i=1}^3\frac{A^2}{(t+1)^2}(\eta_{i}+\tilde{\eta_{i}})^2\\
    +\sum\limits_{j=1}^3\frac{B^2}{(t+1)^2}(\eta_{j}-\tilde{\eta_{j}})^2
\end{multline*}
to the conical metric
$$ds^2+\frac{s^2}{36}\sum\limits_{i=1}^3(\eta_{i}+\tilde{\eta_{i}})^2+
\frac{s^2}{12}\sum\limits_{j=1}^3(\eta_{j}-\tilde{\eta_{j}})^2.$$

Because the metrics are expressed in the same basis to finish the
proof it's sufficient to keep in mind \ref{convergence1} and
remember that $A=2BB_{x}.$ $\Box$

\textbf{Remark}: \textit{Conditions of the theorem \ref{theorem}
could be rewritten as some conditions on initial $G_{2}$-structure
(or metric) for the flow $$\frac{\partial\varphi}{\partial t}\wedge
dr=d\varphi,$$ so we can establish Cauchy problem for the system
\begin{equation}
    \begin{cases}
        \dot B+B^{\prime}=\frac{A}{B}
        \\
        \dot A+A^{\prime}=\frac{1}{2}(1-\frac{A^2}{B^2})
    \end{cases}
\end{equation}}

\newpage

\section*{Appendix}

Let's compute the De Rham differential $d$ of $G_{2}$-structure
$\varphi$ on $C(S^3\times S^3$). To simplify notations we will
denote by $B\prime$ and $A\prime$ partial derivatives
$\frac{\partial B}{\partial r}$ and $\frac{\partial A}{\partial r}$
correspondingly. Recall that in coordinates described in the section
\ref{section}
$$\varphi=e^{456}+e^{621}+e^{174}+e^{527}+e^{637}+e^{135}+e^{432},$$
$$d\varphi=de^{456}+de^{621}+de^{174}+de^{527}+de^{637}+de^{135}+de^{432}.$$

$de^{456}=de^4\wedge e^5\wedge e^6
          -e^4\wedge de^5\wedge e^6
          +e^4\wedge e^5\wedge de^6=
(\frac{B\prime}{B}e^7\wedge e^4+2B(-\eta_{2}\wedge\eta_{3}+
\tilde{\eta_2}\wedge\tilde{\eta_3}))\wedge e^5\wedge e^6-
e^4\wedge(\frac{B\prime}{B}e^7\wedge e^5+
2B(-\eta_3\wedge\eta_1+\tilde{\eta_3}\wedge\tilde{\eta_1}))\wedge
e^6+ e^4\wedge e^5\wedge(\frac{B\prime}{B}e^7\wedge e^6+
2B(-\eta_1\wedge\eta_2+\tilde{\eta_1}\wedge\tilde{\eta_2}))=
3\frac{B\prime}{B}e^7\wedge e^4\wedge e^5\wedge e^6+
\frac{B}{2}(-(\frac{e^2}{A}+\frac{e^5}{B})\wedge(\frac{e^3}{A}+\frac{e^6}{B})
+(\frac{e^2}{A}-\frac{e^5}{B})\wedge(\frac{e^3}{A}-\frac{e^6}{B}))\wedge
e^5\wedge e^6
-\frac{B}{2}e^4\wedge(-(\frac{e^3}{A}+\frac{e^6}{B})\wedge(\frac{e^1}{A}+\frac{e^4}{B})
+(\frac{e^3}{A}-\frac{e^6}{B})\wedge(\frac{e^1}{A}-\frac{e^4}{B}))\wedge
e^6+\frac{B}{2}e^4\wedge e^5\wedge
(-(\frac{e^1}{A}+\frac{e^4}{B})\wedge(\frac{e^2}{A}+\frac{e^5}{B})
+(\frac{e^1}{A}-\frac{e^4}{B})\wedge(\frac{e^2}{A}-\frac{e^5}{B}))=-3\frac{B\prime}{B}e^4\wedge
e^5\wedge e^6\wedge e^7$

$de^{621}=de^6\wedge e^2\wedge e^1
          -e^6\wedge de^2\wedge e^1
          +e^6\wedge e^2\wedge de^1=
(\frac{B\prime}{B}e^7\wedge e^6+2B(-\eta_{1}\wedge\eta_{2}+
\tilde{\eta_1}\wedge\tilde{\eta_2}))\wedge e^2\wedge e^1-
e^6\wedge(\frac{A\prime}{A}e^7\wedge e^2+
2A(-\eta_3\wedge\eta_1-\tilde{\eta_3}\wedge\tilde{\eta_1}))\wedge
e^1+ e^6\wedge e^2\wedge(\frac{A\prime}{A}e^7\wedge e^1+
2A(-\eta_2\wedge\eta_3-\tilde{\eta_2}\wedge\tilde{\eta_3}))=
\frac{B\prime}{B}e^7\wedge e^6\wedge e^2\wedge
e^1+2\frac{A\prime}{A}e^7\wedge e^6\wedge e^2\wedge e^1+
\frac{B}{2}(-(\frac{e^1}{A}+\frac{e^4}{B})\wedge(\frac{e^2}{A}+\frac{e^5}{B})+
(\frac{e^1}{A}-\frac{e^4}{B})\wedge(\frac{e^2}{A}-\frac{e^5}{B}))\wedge
e^2\wedge e^1
-\frac{A}{2}e^6\wedge(-(\frac{e^3}{A}+\frac{e^6}{B})\wedge(\frac{e^1}{A}+\frac{e^4}{B})
-(\frac{e^3}{A}-\frac{e^6}{B})\wedge(\frac{e^1}{A}-\frac{e^4}{B}))\wedge
e^1+\frac{A}{2}e^6\wedge e^2\wedge
(-(\frac{e^2}{A}+\frac{e^5}{B})\wedge(\frac{e^3}{A}+\frac{e^6}{B})
-(\frac{e^2}{A}-\frac{e^5}{B})\wedge(\frac{e^3}{A}-\frac{e^6}{B}))=\frac{B\prime}{B}e^1\wedge
e^2\wedge e^6\wedge e^7+2\frac{A\prime}{A}e^1\wedge e^2\wedge
e^6\wedge e^7$

$de^{174}=de^1\wedge e^7\wedge e^4
          -e^1\wedge de^7\wedge e^4
          +e^1\wedge e^7\wedge de^4=
(\frac{A\prime}{A}e^7\wedge e^1+2A(-\eta_{2}\wedge\eta_{3}-
\tilde{\eta_2}\wedge\tilde{\eta_3}))\wedge e^7\wedge e^4+ e^1\wedge
e^7\wedge(\frac{B\prime}{B}e^7\wedge e^4+
2B(-\eta_2\wedge\eta_3+\tilde{\eta_2}\wedge\tilde{\eta_3}))=
\frac{A}{2}(-(\frac{e^2}{A}+\frac{e^5}{B})\wedge(\frac{e^3}{A}+\frac{e^6}{B})-
(\frac{e^2}{A}-\frac{e^5}{B})\wedge(\frac{e^3}{A}-\frac{e^6}{B}))\wedge
e^7\wedge e^4 +\frac{B}{2}e^1\wedge
e^7\wedge(-(\frac{e^2}{A}+\frac{e^5}{B})\wedge(\frac{e^3}{A}+\frac{e^6}{B})
+(\frac{e^2}{A}-\frac{e^5}{B})\wedge(\frac{e^3}{A}-\frac{e^6}{B}))=-\frac{1}{A}e^2\wedge
e^3\wedge e^7\wedge e^4-\frac{A}{B^2}e^5\wedge e^6\wedge e^7\wedge
e^4-\frac{1}{A}e^1\wedge e^7\wedge e^2\wedge e^6-
\frac{1}{A}e^1\wedge e^7\wedge e^5\wedge e^3= \frac{1}{A}e^2\wedge
e^3\wedge e^4\wedge e^7+\frac{A}{B^2}e^4\wedge e^5\wedge e^6\wedge
e^7-\frac{1}{A}e^1\wedge e^2\wedge e^6\wedge e^7+
\frac{1}{A}e^1\wedge e^3\wedge e^5\wedge e^7$

$de^{527}=de^5\wedge e^2\wedge e^7
          -e^5\wedge de^2\wedge e^7
          +e^5\wedge e^2\wedge de^7=
(\frac{B\prime}{B}e^7\wedge e^5+2B(-\eta_{3}\wedge\eta_{1}+
\tilde{\eta_3}\wedge\tilde{\eta_1}))\wedge e^2\wedge e^7- e^5\wedge
(\frac{A\prime}{A}e^7\wedge e^2+
2A(-\eta_3\wedge\eta_1-\tilde{\eta_3}\wedge\tilde{\eta_1}))\wedge
e^7=
\frac{B}{2}(-(\frac{e^3}{A}+\frac{e^6}{B})\wedge(\frac{e^1}{A}+\frac{e^4}{B})+
(\frac{e^3}{A}-\frac{e^6}{B})\wedge(\frac{e^1}{A}-\frac{e^4}{B}))\wedge
e^2\wedge e^7
+\frac{A}{2}e^5\wedge((\frac{e^3}{A}+\frac{e^6}{B})\wedge(\frac{e^1}{A}+\frac{e^4}{B})
+(\frac{e^3}{A}-\frac{e^6}{B})\wedge(\frac{e^1}{A}-\frac{e^4}{B}))\wedge
e^7=-\frac{1}{A}e^3\wedge e^4\wedge e^2\wedge
e^7-\frac{1}{A}e^6\wedge e^1\wedge e^2\wedge
e^7+\frac{1}{A}e^5\wedge e^3\wedge e^1\wedge e^7+
\frac{A}{B^2}e^5\wedge e^6\wedge e^4\wedge e^7=-\frac{1}{A}e^2\wedge
e^3\wedge e^4\wedge e^7-\frac{1}{A}e^1\wedge e^2\wedge e^6\wedge
e^7-\frac{1}{A}e^1\wedge e^3\wedge e^5\wedge e^7+
\frac{A}{B^2}e^4\wedge e^5\wedge e^6\wedge e^7$

$de^{637}=de^6\wedge e^3\wedge e^7
          -e^6\wedge de^3\wedge e^7
          +e^6\wedge e^3\wedge de^7=
(\frac{B\prime}{B}e^7\wedge e^6+2B(-\eta_{1}\wedge\eta_{2}+
\tilde{\eta_1}\wedge\tilde{\eta_2}))\wedge e^3\wedge e^7- e^6\wedge
(\frac{A\prime}{A}e^7\wedge e^3+
2A(-\eta_1\wedge\eta_2-\tilde{\eta_1}\wedge\tilde{\eta_2}))\wedge
e^7=
\frac{B}{2}(-(\frac{e^1}{A}+\frac{e^4}{B})\wedge(\frac{e^2}{A}+\frac{e^5}{B})+
(\frac{e^1}{A}-\frac{e^4}{B})\wedge(\frac{e^2}{A}-\frac{e^5}{B}))\wedge
e^3\wedge e^7
+\frac{A}{2}e^6\wedge((\frac{e^1}{A}+\frac{e^4}{B})\wedge(\frac{e^2}{A}+\frac{e^5}{B})
+(\frac{e^1}{A}-\frac{e^4}{B})\wedge(\frac{e^2}{A}-\frac{e^5}{B}))\wedge
e^7=-\frac{1}{A}e^1\wedge e^5\wedge e^3\wedge
e^7-\frac{1}{A}e^4\wedge e^2\wedge e^3\wedge
e^7+\frac{1}{A}e^6\wedge e^1\wedge e^2\wedge e^7+
\frac{A}{B^2}e^6\wedge e^4\wedge e^5\wedge e^7=\frac{1}{A}e^1\wedge
e^3\wedge e^5\wedge e^7-\frac{1}{A}e^2\wedge e^3\wedge e^4\wedge
e^7+\frac{1}{A}e^1\wedge e^2\wedge e^6\wedge e^7+
\frac{A}{B^2}e^4\wedge e^5\wedge e^6\wedge e^7$

$de^{135}=de^1\wedge e^3\wedge e^5
          -e^1\wedge de^3\wedge e^5
          +e^1\wedge e^3\wedge de^5=
(\frac{A\prime}{A}e^7\wedge e^1+2A(-\eta_{2}\wedge\eta_{3}-
\tilde{\eta_2}\wedge\tilde{\eta_3}))\wedge e^3\wedge e^5-
e^1\wedge(\frac{A\prime}{A}e^7\wedge e^3+
2A(-\eta_1\wedge\eta_2-\tilde{\eta_1}\wedge\tilde{\eta_2}))\wedge
e^5+ e^1\wedge e^3\wedge(\frac{B\prime}{B}e^7\wedge e^5+
2B(-\eta_3\wedge\eta_1+\tilde{\eta_3}\wedge\tilde{\eta_1}))=
\frac{B\prime}{B}e^7\wedge e^1\wedge e^3\wedge
e^5+2\frac{A\prime}{A}e^7\wedge e^1\wedge e^3\wedge e^5
+\frac{A}{2}(-(\frac{e^2}{A}+\frac{e^5}{B})\wedge(\frac{e^3}{A}+\frac{e^6}{B})-
(\frac{e^2}{A}-\frac{e^5}{B})\wedge(\frac{e^3}{A}-\frac{e^6}{B}))\wedge
e^3\wedge e^5
+\frac{A}{2}e^1\wedge((\frac{e^1}{A}+\frac{e^4}{B})\wedge(\frac{e^2}{A}+\frac{e^5}{B})
+(\frac{e^1}{A}-\frac{e^4}{B})(\frac{e^2}{A}-\frac{e^5}{B}))\wedge
e^5+\frac{B}{2}e^1\wedge e^3\wedge
(-(\frac{e^3}{A}+\frac{e^6}{B})\wedge(\frac{e^1}{A}+\frac{e^4}{B})
+(\frac{e^3}{A}-\frac{e^6}{B})\wedge(\frac{e^1}{A}-\frac{e^4}{B}))=-\frac{B\prime}{B}e^1\wedge
e^3\wedge e^5\wedge e^7-2\frac{A\prime}{A}e^1\wedge e^3\wedge
e^5\wedge e^7$

$de^{432}=de^4\wedge e^3\wedge e^2
          -e^4\wedge de^3\wedge e^2
          +e^4\wedge e^3\wedge de^2=
(\frac{B\prime}{B}e^7\wedge e^4+2B(-\eta_{2}\wedge\eta_{3}+
\tilde{\eta_2}\wedge\tilde{\eta_3}))\wedge e^3\wedge e^2-
e^4\wedge(\frac{A\prime}{A}e^7\wedge e^3+
2A(-\eta_1\wedge\eta_2-\tilde{\eta_1}\wedge\tilde{\eta_2}))\wedge
e^2+ e^4\wedge e^3\wedge(\frac{A\prime}{A}e^7\wedge e^2+
2A(-\eta_3\wedge\eta_1-\tilde{\eta_3}\wedge\tilde{\eta_1}))=
\frac{B\prime}{B}e^7\wedge e^4\wedge e^3\wedge
e^2+2\frac{A\prime}{A}e^7\wedge e^4\wedge e^3\wedge e^2
+\frac{B}{2}(-(\frac{e^2}{A}+\frac{e^5}{B})\wedge(\frac{e^3}{A}+\frac{e^6}{B})
+(\frac{e^2}{A}-\frac{e^5}{B})\wedge(\frac{e^3}{A}-\frac{e^6}{B}))\wedge
e^3\wedge e^2
+\frac{A}{2}e^4\wedge((\frac{e^1}{A}+\frac{e^4}{B})\wedge(\frac{e^2}{A}+\frac{e^5}{B})
+(\frac{e^1}{A}-\frac{e^4}{B})\wedge(\frac{e^2}{A}-\frac{e^5}{B}))\wedge
e^2-\frac{A}{2}e^4\wedge e^3\wedge
((\frac{e^3}{A}+\frac{e^6}{B})\wedge(\frac{e^1}{A}+\frac{e^4}{B})
+(\frac{e^3}{A}-\frac{e^6}{B})\wedge(\frac{e^1}{A}-\frac{e^4}{B}))=\frac{B\prime}{B}e^2\wedge
e^3\wedge e^4\wedge e^7+2\frac{A\prime}{A}e^2\wedge e^3\wedge
e^4\wedge e^7$

$d\varphi=-3\frac{B\prime}{B}e^4\wedge e^5\wedge e^6\wedge
e^7+\frac{B\prime}{B}e^1\wedge e^2\wedge e^6\wedge
e^7+2\frac{A\prime}{A}e^1\wedge e^2\wedge e^6\wedge
e^7+\frac{1}{A}e^2\wedge e^3\wedge e^4\wedge
e^7+\frac{A}{B^2}e^4\wedge e^5\wedge e^6\wedge
e^7-\frac{1}{A}e^1\wedge e^2\wedge e^6\wedge e^7+
\frac{1}{A}e^1\wedge e^3\wedge e^5\wedge e^7-\frac{1}{A}e^2\wedge
e^3\wedge e^4\wedge e^7-\frac{1}{A}e^1\wedge e^2\wedge e^6\wedge
e^7-\frac{1}{A}e^1\wedge e^3\wedge e^5\wedge e^7+
\frac{A}{B^2}e^4\wedge e^5\wedge e^6\wedge e^7+\frac{1}{A}e^1\wedge
e^3\wedge e^5\wedge e^7-\frac{1}{A}e^2\wedge e^3\wedge e^4\wedge
e^7+\frac{1}{A}e^1\wedge e^2\wedge e^6\wedge e^7+
\frac{A}{B^2}e^4\wedge e^5\wedge e^6\wedge
e^7-\frac{B\prime}{B}e^1\wedge e^3\wedge e^5\wedge
e^7-2\frac{A\prime}{A}e^1\wedge e^3\wedge e^5\wedge
e^7+\frac{B\prime}{B}e^2\wedge e^3\wedge e^4\wedge
e^7+2\frac{A\prime}{A}e^2\wedge e^3\wedge e^4\wedge
e^7=3(\frac{A}{B^2}-\frac{B\prime}{B})e^4\wedge e^5\wedge e^6\wedge
e^7+(\frac{B\prime}{B}+2\frac{A\prime}{A}-\frac{1}{A})e^1\wedge
e^2\wedge e^6\wedge
e^7+(\frac{B\prime}{B}+2\frac{A\prime}{A}-\frac{1}{A})e^2\wedge
e^3\wedge e^4\wedge
e^7-(\frac{B\prime}{B}+2\frac{A\prime}{A}-\frac{1}{A})e^1\wedge
e^3\wedge e^5\wedge e^7$
\\Let's compute now $\frac{\partial \varphi}{\partial t}$. Recall that\\
$e^i=A(\eta_{i}+\tilde{\eta_{i}})$ for $i=1,2,3$\\
$e^j=B(\eta_{j}-\tilde{\eta_{j}})$ for $j=4,5,6$\\
$e^7=dr.$\\
Then
\begin{multline*}
$$\\
\frac{\partial e^{456}}{\partial t}=3\frac{\dot B}{B}e^{456}\\
\frac{\partial e^{621}}{\partial
t}=(\frac{\dot B}{B}+2\frac{\dot A}{A})e^{621}\\
\frac{\partial e^{174}}{\partial
t}=(\frac{\dot B}{B}+\frac{\dot A}{A})e^{174}\\
\frac{\partial e^{527}}{\partial
t}=(\frac{\dot B}{B}+\frac{\dot A}{A})e^{527}\\
\frac{\partial e^{637}}{\partial
t}=(\frac{\dot B}{B}+\frac{\dot A}{A})e^{637}\\
\frac{\partial e^{135}}{\partial t}=(\frac{\dot B}{B}+2\frac{\dot A}{A})e^{135}\\
\frac{\partial e^{432}}{\partial t}=(\frac{\dot B}{B}+2\frac{\dot
A}{A})e^{432}\\$$
\end{multline*}
where $\dot A$ and $\dot B$ are partial w.r.t. $t$ derivatives of
functions $A$ and $B$ correspondingly.

\textit{Proof of the lemma \ref{lemma1}:}

Keeping in mind above calculations we have\\
$  3\frac{\dot B}{B}e^{4567}-(\frac{\dot B}{B}+2\frac{\dot
  A}{A})e^{1267}+ (\frac{\dot B}{B}+2\frac{\dot
  A}{A})e^{1357}-(\frac{\dot B}{B}+2\frac{\dot A}{A})e^{2347}
  =\frac{\partial \varphi}{\partial t}\wedge
  e^7=d\varphi=3(\frac{A}{B^2}-\frac{B\prime}{B})e^{4567}+
  (\frac{B\prime}{B}+2\frac{A\prime}{A}-\frac{1}{A})e^{1267}+
  (\frac{B\prime}{B}+2\frac{A\prime}{A}-\frac{1}{A})e^{2347}-
  (\frac{B\prime}{B}+2\frac{A\prime}{A}-\frac{1}{A})e^{1357}
$\\
or\\
\begin{equation*}
    \begin{cases}
        \frac{\dot B}{B}=-\frac{B\prime}{B}+\frac{A}{B^2}
        \\
        2\frac{\dot A}{A}+\frac{\dot B}{B}=\frac{1}{A}-2\frac{A\prime}{A}-\frac{B\prime}{B}
    \end{cases}
\end{equation*}\\
or\\
\begin{equation}\label{system1}
    \begin{cases}
        \dot B+B\prime=\frac{A}{B}
        \\
        \dot A+A\prime=\frac{1}{2}(1-\frac{A^2}{B^2})
    \end{cases}
\end{equation}\\
Let's change variables $t$ and $r$ on $x=r+t$ and $y=r-t$. In these
variables equations \ref{system1} look as follow
\begin{equation}\label{system2}
    \begin{cases}
        2B_{x}=\frac{A}{B}
        \\
        2A_{x}=\frac{1}{2}(1-\frac{A^2}{B^2})
    \end{cases}
\end{equation}\\
By expressing the function $A$ from the first equation and
substituting it into the second equation of \ref{system2} we get the
desired system
\begin{equation*}
    \begin{cases}
        A=2BB_{x}
        \\
        8BB_{xx}+12B_{x}^2=1
    \end{cases}
\end{equation*} $\Box$

\newpage
\bibliography{scibib}
\bibliographystyle{Science}

\textsc{Sobolev Institute of Mathematics, Novosibirsk, Russia}

\textbf{hazhgaly@gmail.com}

\end{document}